\documentclass{article}
\usepackage{amsmath,amsfonts,amsthm,amssymb,amscd}

\newcommand{\esssup}{\mathop{\rm ess\ sup}}

\newcommand{\p}{\partial}
\newcommand{\e}{\varepsilon}

\newcommand{\R}{{\mathbb R}}

\newcommand{\T}{{\mathbb T}}

\newcommand{\nnn}{{\boldsymbol{\mathit n}}}

\newcommand{\CC}{{\cal C}}
\newcommand{\DD}{{\cal D}}
\newcommand{\EE}{{\cal E}}
\newcommand{\FF}{{\cal F}}

\newcommand{\KK}{{\cal K}}

\newcommand{\RR}{{\cal R}}
\newcommand{\sS}{{\cal S}}

\newcommand{\XX}{{\cal X}}

\newcommand{\diver}{\mathop{\rm div}\nolimits}

\newcommand{\dist}{\mathop{\rm dist}\nolimits}

\theoremstyle{plain}
\newtheorem{theorem}{Theorem}[section]

\newtheorem{lemma}[theorem]{Lemma}
\newtheorem{proposition}[theorem]{Proposition}
\newtheorem{corollary}[theorem]{Corollary}
\theoremstyle{definition}
\newtheorem{definition}[theorem]{Definition}

\theoremstyle{remark}

\numberwithin{equation}{section}

\begin{document}
\author{Armen Shirikyan} 
\title{Exact controllability in projections 
for three-dimensional Navier--Stokes equations}
\date{\small Laboratoire de Math\'ematiques\\
Universit\'e de Paris-Sud XI, B\^atiment 425\\ 
91405 Orsay Cedex, France\\ 
E-mail: Armen.Shirikyan@math.u-psud.fr}
\maketitle
\begin{abstract}
The paper is devoted to studying controllability properties for 3D
Navier--Stokes equations in a bounded domain. We establish a
sufficient condition under which the problem in question is exactly
controllable in any finite-dimensional projection. Our sufficient
condition is verified for any torus in~$\R^3$. The proofs are based on
a development of a general approach introduced by Agrachev and
Sarychev in the 2D case. As a simple consequence of the result on
controllability, we show that the Cauchy problem for the 3D
Navier-Stokes system has a unique strong solution for any initial
function and a large class of external forces.


\end{abstract}
\tableofcontents
\setcounter{section}{-1}

\section{Introduction}
\label{s0}
Let us consider the three-dimensional Navier--Stokes (NS) system 
\begin{equation} \label{1}
\dot u+(u,\nabla)u-\nu\Delta u+\nabla p=f(t,x),\quad \diver u=0,
\end{equation}
where the space variables $x=(x_1,x_2,x_3)$ belong to a
three-dimensional torus $\T^3\subset\R^3$, $\nu>0$~is the viscosity,
$u=(u_1,u_2,u_3)$ and~$p$ are unknown velocity field and pressure,
and~$f(t,x)$ is an external force. Suppose that~$f$ is represented in
the form
\begin{equation} \label{2}
f(t,x)=h(t,x)+\eta(t,x),
\end{equation}
where $h$ is a given function and~$\eta$ is a control taking on values
in a {\it finite-dimensional\/} subspace $E\subset L^2(\T^3,\R^3)$. 
Equations~\eqref{1}, \eqref{2} are supplemented
with the initial condition
\begin{equation} \label{3}
u(0)=u_0,
\end{equation}
where $u_0\in H^1(\T^3,\R^3)$ is a divergence-free vector field. Let
us denote by~$H$ the space of functions $u\in L^2(\T^3,\R^3)$ such
that $\diver u=0$ on~$\T^3$. We fix an arbitrary subspace~$F\subset H$
and denote by ${\mathsf P}_F:H\to H$ the orthogonal projection
onto~$F$. Problem~\eqref{1}, \eqref{2} is said to be {\it controllable
  in a time~$T>0$ for the projection to~$F$} if for any initial
function~$u_0$ and any~$\hat u\in F$ there exists an infinitely smooth
control $\eta:[0,T]\to E$ such that~\eqref{1}~--~\eqref{3} has a
unique strong solution~$u(t;\eta)$, which satisfies the relation
\begin{equation} \label{4}
{\mathsf P}_Fu(T;\eta)=\hat u.
\end{equation}
One of the main results of this paper says that if the space~$E$ is
sufficiently large, then for any~$T>0$ and~$\nu>0$ and any
finite-dimensional subspace~$F\subset H$ problem~\eqref{1}, \eqref{2}
is controllable in time~$T$ for the projection to~$F$.

A general approach for studying controllability of PDE's in
finite-dimensional projections was introduced by Agrachev and Sarychev
in the landmark article~\cite{AS-2004} (see also~\cite{AS-2006}). They
considered the 2D NS system on a torus and proved that it is
controllable for the projection to any finite-dimensional space~$F$,
with a control function taking on values in a fixed subspace~$E$. We
emphasise that the time of control~$T$ can be chosen arbitrarily
small, and the control space~$E$ does not depend on~$\nu$
and~$T$.\,\footnote{It is shown in~\cite{AS-2006} that the 2D
Euler and Navier--Stokes equations are controllable by a control of
dimension four.} The Agrachev--Sarychev approach is based on the
concept of {\it solid controllability\/} (cf.~Definition~\ref{d2.6} of
the present paper). They construct explicitly an increasing sequence
of finite-dimensional subspaces~$\{E_k\}_{k\ge0}$ such that $E_0=E$,
and the following two properties hold.
\begin{itemize}
\item[\bf(i)]
There is an integer $N\ge1$ such that the NS system is solidly
controllable by an $E_N$-valued control.
\item[\bf(ii)]
If the NS system is solidly controllable by an $E_k$-valued control for
some integer~$k\ge1$, then it is solidly controllable by an $E_{k-1}$-valued
control. 
\end{itemize}
These two assertions imply the required result. 

\smallskip
In this paper, we take a slightly different viewpoint based on uniform
approximate controllability.\,\footnote{Note that the concept of
uniform approximate controllability is implicitly present in the
Agrachev--Sarychev argument~\cite{AS-2004}.} Namely, we shall say that
the NS system~\eqref{1}, \eqref{2} is {\it uniformly approximately
controllable\/} (UAC) if for any constant~$\e>0$, any initial
function~$u_0$, and any compact subset~$\KK$ of the phase space there
is a continuous mapping~$\varPsi$ from~$\KK$ to the space of
$E$-valued controls such that for every~$\hat u\in\KK$
problem~\eqref{1}~--~\eqref{3} with $\eta=\varPsi(\hat u)$ has a
unique strong solution~$u(t;\eta)$, which satisfies the inequality
\begin{equation} \label{5}
\|u(T;\eta)-\hat u\|<\e,
\end{equation}
where $\|\cdot\|$ denotes the $L^2$-norm. It turns out that
assertions~(i) and~(ii) remain valid for the 3D NS system if we
replace the solid controllability by uniform approximate
controllability (cf.~\cite{S05}). Hence, we prove that if~$E$ is
sufficiently large, then problem~\eqref{1}, \eqref{2} is UAC by an
$E$-valued control. The required result on exact controllability in
finite-dimensional projections is a simple consequence of the above
property. Indeed, let~$B_F(R)$ be the closed ball in~$F$ of radius~$R$
centred at origin and let $\KK=B_F(R)$. In this case, it follows
from~\eqref{5} that 
\begin{equation} \label{6}
\|{\mathsf P}_Fu(T;\varPsi(\hat u))-\hat u\|<\e\quad
\mbox{for any $\hat u\in B_F(R)$}.
\end{equation}
The function $\varPhi:\hat u\mapsto{\mathsf P}_Fu(T,\varPsi(\hat u))$
is continuous from~$B_F(R)$ to~$F$. Using the Brouwer fixed point
theorem and inequality~\eqref{6}, it is easy to show (see
Proposition~\ref{p1.1}) that $\varPhi(B_F(R))$ contains the
ball~$B_F(R-\e)$. Since~$R>0$ is arbitrary, we conclude that~\eqref{4}
holds for any~$\hat u\in F$ and an appropriate $E$-valued control
function~$\eta$.

In conclusion, we note that the problem of controllability and
stabilisation for the Navier--Stokes and Euler equation was in the
focus of attention of many researchers; for instance, see the
papers~\cite{fursikov-1995,coron-1995,coron-1996,CF-1996,imanuilov-1998,
FE-1999,coron-1999,FC-1999,glass-2000,BS-2001,fursikov-2001,fursikov-2004,BT-2004}
and references therein. However, the powerful techniques developed in
those papers do not apply to the present setting because of the
specific type of control we are interested in.

\smallskip
The paper is organised as follows. In Section~\ref{s1}, we recall a
simple sufficient condition for surjectivity of continuous mapping in
a finite-dimensional space and formulate two perturbative results on
unique solvability of NS-type equations. Section~\ref{s2} contains the
formulations of the main results of this paper. We also discuss some
corollaries on solid controllability in finite-dimensional projection
and the Cauchy problem for the 3D NS system. The proofs are presented
in Section~\ref{s3}. Finally, in the Appendix, we prove an auxiliary
result used in Section~\ref{s3}.

\medskip
{\bf Acknowledgements}. I am grateful to A.~A.~Agrachev for
stimulating discussion on the subject. 

\medskip
{\bf Notation}.  We denote by~$\R_+$ the half-line $[0,+\infty)$ and
by~$J_T$ the interval~$[0,T]$. If~$s\ge1$ and~$r\ge0$ are some
integers, then we set $J_T(r,s)=[t_r,t_{r+1})$, where
$t_r=rT/s$. Let $J\subset\R_+$ be a closed interval, let
$D\subset\R^3$ be a bounded domain, let~$X$ be a Banach space with a
norm~$\|\cdot\|$, and let~$\KK$ be a metric space. We shall use the
following functional spaces.

\medskip
\noindent
$L^p(J,X)$ is the space of measurable functions $u:J\to X$
with finite norm
\begin{equation} \label{11}
\|u\|_{L^p(J,X)}:=\biggl(\int_J\|u(t)\|_X^pdt\biggr)^{1/p},
\end{equation}
where $\|\cdot\|_X$ stands for the norm in~$X$. If $p=\infty$,
then~\eqref{11} is replaced by
$$
\|u\|_{L^\infty(J,X)}:=\esssup_{t\in J}\|u(t)\|_X.
$$
We shall write $L^p(J)$ instead of~$L^p(J,\R)$.

\smallskip
\noindent
$L_{\rm loc}^p(\R_+,X)$ is the space of functions $u:R_+\to X$ whose
restriction to any interval $J\subset\R_+$ belongs to~$L^p(J,X)$. 

\smallskip
\noindent
$C^k(J,X)$ is the space of continuous functions $u:J\to X$ that
are~$k$ times continuously differentiable. In the case $k=0$, we shall
write~$C(J,X)$.

\smallskip
\noindent
$C(\KK,X)$ is the space of continuous functions $u:\KK\to X$. If
$X=\R$, then we write~$C(\KK)$.

\smallskip
\noindent
$H^s(D,\R^3)$ is the space of vector functions $(u_1,u_2,u_3)$ whose
components belong to the Sobolev space of order~$s$. In the case $s=0$,
it coincides with the Lebesgue space~$L^2(D,\R^3)$.  

\smallskip
\noindent
$H$, $V$, $U$, and~$\XX_T$ are standard functional spaces arising in
the theory of Navier-Stokes equations; they are defined in
Section~\ref{s1.2}.

\section{Preliminaries}
\label{s1}
%We have compiled here some auiximiary results that are
%used in Sections~\ref{s2} and~\ref{s3}. 

\subsection{Image of continuous mappings}
\label{s1.1}
Let~$F$ be a finite-dimensional vector space  with a norm~$\|\cdot\|_F$, 
let~$B_F(R)$ be the closed ball in~$F$ of radius~$R$ centred at origin, 
and let $\varPhi:B_F(R)\to F$ be a continuous mapping. The following 
result is a simple consequence of the Brouwer theorem. 

\begin{proposition} \label{p1.1}
Suppose there is a constant $\e\in(0,R)$ such that
\begin{equation} \label{1.1}
\|\varPhi(u)-u\|_F\le\e\quad\mbox{for any $u\in B_F(R)$}.
\end{equation}
Then $\varPhi(B_F(R))\supset B_F(R-\e)$.
\end{proposition}

\begin{proof}
Let us fix any point $\hat u\in B_F(R-\e)$ and consider the continuous
mapping 
$$
\varPsi:B_F(R)\to F, \quad\varPsi(u)=\hat u-\varPhi(u)+u.
$$
It follows from~\eqref{1.1} that $\varPsi(B_F(R))\subset B_F(R)$. 
Therefore, by the Brouwer theorem (e.g., see Section~1.19 in~\cite{Taylor}),
$\varPsi$~has a fixed point $u_0\in B_F(R)$. Direct verification shows 
that $\varPhi(u_0)=\hat u$. Thus, any point $\hat u\in B_F(R-\e)$ has a 
preimage, and we obtain the required inclusion.
\end{proof}

\subsection{Strong solutions of Navier--Stokes type equations}
\label{s1.2}
We first introduce some standard functional spaces arising in the
theory of 3D Navier--Stokes (NS) equations. Let
$$
H=\bigl\{u\in L^2(D,\R^3):
\mbox{$\diver u=0$ in~$D$, $(u,\nnn)|_{\p D}=0$}\bigr\},
$$
where~$\nnn$ is the outward unit normal to~$\p D$, and let~$\Pi$ be
the orthogonal projection in~$L^2(D,\R^3)$ onto the closed
subspace~$H$. We denote by~$H^s=H^s(D,\R^3)$ the space of vector
functions $u=(u_1,u_2,u_3)$ with components in the Sobolev class of
order~$s$ and by~$H^s_0(D,\R^3)$ the space of functions $u\in H^s$
vanishing on~$\p D$. Let~$\|\cdot\|_s$ be the usual norm
in~$H^s$. In the case $s=0$, we write~$\|\cdot\|$. Define the spaces
$$
V=H_0^1(D,\R^3)\cap H,\quad U=H^2(D,\R^2)\cap V
$$
and endow them with natural norms. 

It is well known (e.g., see~\cite{Te}) that the NS system is
equivalent to the following evolution equation in~$H$:  
\begin{equation} \label{1.2}
\dot u+\nu Lu+B(u)=f(t), 
\end{equation}
where $L=-\Pi\Delta$ is the Stokes operator and
$B(u)=\Pi\{(u,\nabla)u\}$ is the bilinear form resulting from the
nonlinear term in the original system. Let $E\subset U$ be a
finite-dimensional vector space and let~$E^\bot$ be its orthogonal
complement in~$H$. Denote by ${\mathsf P}={\mathsf P}_E$ and
${\mathsf Q}={\mathsf Q}_E$ the orthogonal projections in~$H$ onto the
subspaces~$E$ and~$E^\bot$, respectively. Along with~\eqref{1.2},
consider the Cauchy problem
\begin{align}
\dot w+\nu L_Ew+{\mathsf Q}\bigl(B(w)+B(v,w)+B(w,v)\bigr)&=f(t), 
\label{1.3}\\
w(0)&=w_0, \label{1.4}
\end{align}
where $L_E={\mathsf Q}L$, $B(v,w)=\Pi\{(v,\nabla)w\}$, 
and $v\in L^4(J_T,V)$ and $f\in L^2(J_T,E^\bot)$ are given functions.
We set
$$
\XX_T=C(J_T,V)\cap L^2(J_T,U),\quad 
\XX_T(E)=C(J_T,V\cap E^\bot)\cap L^2(J_T,U\cap E^\bot).
$$
The following result is established in~\cite[Section~1.4]{S05} (see
Theorem~1.8). 

\begin{proposition} \label{p1.2}
For any $\nu>0$ and $R>0$ there are positive constants~$\e$ and~$C$
such that the following assertions hold.
\begin{itemize}
\item[{\bf(i)}]
Let $\hat v\in L^4(J_T,H^1)$, $\hat f\in L^2(J_T,E^\bot)$, and
$\widehat w_0\in V\cap E^\bot$ be some functions such that
problem~\eqref{1.3}, \eqref{1.4} with $v=\hat v$, $f=\hat f$,
$w_0=\widehat w_0$ has a solution $\widehat w\in\XX_T(E)$. Suppose
that
$$ %\begin{equation}  \label{1.5}
\|\hat v\|_{L^4(J_T,H^1)}\le R, \quad \|\hat f\|_{L^2(J_T,E^\bot)}\le R,
\quad \|\widehat w\|_{\XX_T}\le R.
$$ %\end{equation}
Then, for any triple $(v,f,w_0)$ satisfying the inequalities
\begin{equation}  \label{1.6}
\|v-\hat v\|_{L^4(J_T,H^1)}\le\e, \quad \|f-\hat f\|_{L^2(J_T,E^\bot)}\le\e,
\quad \|w_0-\widehat w_0\|_V\le\e,
\end{equation}
problem~\eqref{1.3}, \eqref{1.4} has a unique
solution~$w\in\XX_T(E)$.
\item[{\bf(ii)}]
Let
$$
\RR: L^4(J_T,H^1)\times L^2(J_T,E^\bot)\times(V\cap E^\bot)\to\XX_T(E)
$$
be an operator that is defined on the set of functions~$(v,f,w_0)$
satisfying~\eqref{1.6} and takes each triple $(v,f,w_0)$ to the
solution~$w\in\XX_T(E)$ of~\eqref{1.3}, \eqref{1.4}.  Then~$\RR$ is
uniformly Lipschitz continuous, and its Lipschitz constant does not
exceed~$C$.
\end{itemize}
\end{proposition}

We now consider Eq.~\eqref{1.3}  in which~$E$ is a finite-dimensional
vector space spanned by some eigenfunctions of the Stokes
operator~$L$. Namely, let~$\{e_j\}$ be a complete set of normalised
eigenfunctions for~$L$, let~$H_N$ be the vector span of 
$\{e_j,1\le j\le N\}$, and let~$H_N^\bot$ be the orthogonal complement
of~$H_N$ in the space~$H$. We denote by~${\mathsf P}_N$ 
and~${\mathsf Q}_N$ the orthogonal projections in~$H$ onto the
subspaces~$H_N$ and~$H_N^\bot$, respectively. 

Let us consider the equation
\begin{equation} \label{1.7}
\dot w+\nu L_N(w+v)+{\mathsf Q}_NB(w+v)=f(t),
\end{equation}
where $L_N={\mathsf Q}_NL$.

\begin{proposition} \label{p1.4}
For any $R>0$ and~$\nu>0$ there is an integer $N_0\ge1$ and a
constant~$C>0$ such that the following assertions hold.
\begin{itemize}
\item[\bf(i)]
Let $N\ge N_0$ be an integer and let functions $v\in\XX_T$, 
$f\in L^2(J_T,H_N^\bot)$, and $w_0\in H_N^\bot\cap V$ satisfy the
inequalities
\begin{equation} \label{1.8}
\|v\|_{\XX_T}\le R, \quad \|f\|_{L^2(J_T,H)}\le R, 
\quad \|w_0\|_V\le R.
\end{equation}
Then problem~\eqref{1.7}, \eqref{1.4} has a unique solution
$w\in\XX_T(H_N)$.
\item[\bf(ii)] 
Let~$\sS$ be an operator that takes each triple $(v,f,w_0)$
satisfying~\eqref{1.8} to the solution~$w\in\XX_T(H_N)$ of~\eqref{1.7},
\eqref{1.4}. Then~$\sS$ is uniformly Lipschitz continuous in the
corresponding spaces, and its Lipschitz constant does not exceed~$C$.
\end{itemize}
\end{proposition}

\begin{proof}
The existence and uniqueness of solution is established in~\cite{S05}
(see Proposition~1.10). The proof of~(ii) is rather standard, and we
only outline it. 

Let $w^1,w^2\in\XX_T(H_N)$ be two solutions of~\eqref{1.7},
\eqref{1.4} that correspond to some triples~$(v^i,f^i,w_0^i)$,
$i=1,2$. Then the function $w=w^1-w^2\in\XX_T(H_N)$ is a solutions of
the problem
$$
\dot w+\nu L_Nw=g(t), \quad w(0)=w_0^1-w_0^2,
$$
where we set
$$
g(t)=(f^1-f^2)-\nu L_N(v^1-v^2)-{\mathsf Q}_N(B(w^1+v^1)-B(w^2+v^2)).
$$
Repeating literally the argument used in Step~2 of the proof of
Proposition~1.10 in~\cite{S05}, we show that if~$N$ is sufficiently
large, then
\begin{equation} \label{1.9}
\|w\|_{\XX_T}\le 
C_1\bigl(\|f^1-f^2\|_{L^2(J_T,H)}+\|v^1-v^2\|_{L^2(J_T,U)}
+\|w_0^1-w_0^2\|_V\bigr)+\tfrac12\|w\|_{\XX_T},
\end{equation}
where $C_1>0$ is a constant depending only $T$, $R$,
and~$\nu$. Inequality~\eqref{1.9} implies the required result.
\end{proof}

\section{Main results}
\label{s2}

\subsection{Exact controllability in observed projections}
\label{s2.1}
Consider the controlled Navier--Stokes (NS) equations
\begin{align}
\dot u+\nu Lu+B(u)&=h(t)+\eta(t),\label{2.1}\\
u(0)&=u_0, \label{2.2}
\end{align}
where $h\in L_{\rm loc}^2(\R_+,H)$ and $u_0\in V$ are given functions
and~$\eta$ is a control function with range in a finite-dimensional
vector space $E\subset U$. For any $h\in L^2(J_T,H)$, $u_0\in V$,
and~$T>0$, we denote by~$\Theta_T(h,u_0)$ the set of functions
$\eta\in L^2(J_T,H)$ for which problem~\eqref{2.1}, \eqref{2.2} has a
unique solution $u\in\XX_T$.  It follows from Proposition~\ref{p1.2}
with $E=\{0\}$ and $v\equiv0$ that~$\Theta_T(h,u_0)$ is an open subset
of~$L^2(J_T,H)$.

Let us fix a constant~$T>0$, a finite-dimensional space~$F\subset H$,
and a projection ${\mathsf P}_F:H\to H$ onto~$F$.

\begin{definition} \label{d2.1}
Equation~\eqref{2.1} with $\eta\in L^2(J_T,E)$ is said to be {\it
${\mathsf P}_F$-controllable in time~$T$\/} if for any $u_0\in V$
and~$\hat u\in F$ there is $\eta\in\Theta_T(h,u_0)\cap L^2(J_T,E)$
such that
\begin{equation} \label{2.3}
{\mathsf P}_Fu(T)=\hat u,
\end{equation}
where $u\in\XX_T$ denotes the solution of~\eqref{2.1}, \eqref{2.2}. 
\end{definition}

To formulate the main result of this paper, we introduce some notation. 
For any finite-dimensional subspace $E\subset U$, we denote by~$\FF(E)$ 
the largest vector space~$G\subset U$ such that any element 
$\eta_1\in G$ is representable in the form
$$ %\begin{equation} \label{2.4}
\eta_1=\eta-\sum_{j=1}^k\alpha_jB(\zeta^j),
$$ %\end{equation}
where $\eta,\zeta^1,\dots,\zeta^k\in E$ are some vectors and
$\alpha_1,\dots,\alpha_k$ are non-negative constants.  Since~$B$ is
a quadratic operator, we see that $\FF(E)\subset U$ is a well-defined
finite-dimensional subspace containing~$E$. For a finite-dimensional
subspace $E\subset U$, we set
\begin{equation} \label{2.5}
E_0=E, \quad E_k=\FF(E_{k-1})\quad\mbox{for $k\ge1$}, 
\quad E_\infty=\bigcup_{k=1}^\infty E_k.
\end{equation}
The following theorem is the main result of this paper.

\begin{theorem} \label{t2.2}
Let $h\in L_{\rm loc}^2(\R_+,H)$ and let $E\subset U$ be a
finite-dimensional subspace such that~$E_\infty$ is dense in~$H$. Then
for any~$T>0$, any finite-dimensional subspace~$F\subset H$, and any projection
${\mathsf P}_F:H\to H$ onto~$F$ the Navier--Stokes system~\eqref{2.1}
with $\eta\in L^2(J_T,E)$ is ${\mathsf P}_F$-controllable in
time~$T$. Moreover, the control function~$\eta$ can be chosen from the
space~$C^\infty(J_T,E)$.
\end{theorem}

In the case of a general bounded domain, it is difficult to check
whether~$E_\infty$ is dense in~$H$. However, Theorem~\ref{t2.2}
remains valid for the NS equation~\eqref{2.1} on a 3D torus, and it is
shown in~\cite[Section~2.3]{S05} that\footnote{Recall that~$H_N$
denotes the vector space spanned by the first~$N$ eigenfunctions of
the Stokes operator~$L$.} if~$E\supset H_N$ for a sufficiently
large~$N\ge1$, then~$E_\infty$ contains all the eigenfunctions
of~$L$. Thus, we obtain the following result.

\begin{corollary} \label{c2.3}
Let~$\T^3$ be a torus in~$\R^3$. Then there is an integer~$N\ge1$ such
that if the control space~$E$ contains~$H_N$, then for any
constants~$\nu>0$ and~$T>0$, any function $h\in L_{\rm
loc}^2(\R_+,H)$, any finite-dimensional subspace~$F\subset H$, and any
projection ${\mathsf P}_F:H\to H$ onto~$F$ the Navier--Stokes
system~\eqref{2.1} on~$\T^3$ with~$\eta\in L^2(J_T,E)$ is ${\mathsf
P}_F$-controllable in time~$T$, and the control function~$\eta$ can be
chosen from the space~$C^\infty(J_T,E)$.
\end{corollary}

The proof of Theorem~\ref{t2.2} is based on a property of uniform
approximate controllability for~\eqref{2.1}. That concept is of
independent interest and is discussed in the next subsection.

\subsection{Uniform approximate controllability}
\label{s2.2}

Let us fix any $T>0$ and $h\in L^2(J_T,H)$ and denote by $\RR(u_0,\eta)$
an operator that is defined on the set
$$
D(\RR)=\bigl\{(u_0,\eta)\in V\times L^2(J_T,H): 
\eta\in\Theta_T(h,u_0)\bigr\}
$$
and takes each pair~$(u_0,\eta)\in D(\RR)$ to the solution~$u\in\XX_T$
of problem~\eqref{2.1}, \eqref{2.2}. Proposition~\ref{p1.2} with
$E=\{0\}$ and $v\equiv0$ implies that~$D(\RR)$ is an open subset of 
$V\times L^2(J_T,H)$, and~$\RR$ is locally Lipschitz continuous
on~$D(\RR)$. For any $t\in J_T$, we denote by~$\RR_t(u_0,\eta)$ the
restriction of~$\RR(u_0,\eta)$ to the time~$t$.

Let $X\subset L^2(J_T,H)$ be an arbitrary vector space, not necessarily closed. 
We endow~$X$ with the norm of~$L^2(J_T,H)$.

\begin{definition} \label{d2.4}
Equation~\eqref{2.1}  with~$\eta\in X$ is said to be 
{\it uniformly approximately controllable in time~$T$\/} 
if for any initial point~$u_0\in V$,
any compact set~$\KK\subset V$, and any~$\e>0$ there is a continuous function
$$
\varPsi:\KK\to X\cap\Theta_T(h,u_0)
$$
such that
\begin{equation} \label{2.6}
\|\RR_T(u_0,\varPsi(\hat u))-\hat u\|_V<\e\quad
\mbox{for any $\hat u\in\KK$}.
\end{equation}
\end{definition}

The following result shows that, under the conditions of
Theorem~\ref{t2.2}, Eq.~\eqref{2.1} is uniformly approximately
controllable (UAC).

\begin{theorem} \label{t2.5}
Let $h\in L_{\rm loc}^2(\R_+,H)$ and let $E\subset U$ be a
finite-dimensional subspace such that~$E_\infty$ is dense in~$H$. Then
for any~$T>0$ and~$\nu>0$ the Navier--Stokes system~\eqref{2.1} with
$\eta\in C^\infty(J_T,E)$ is UAC in time~$T$.
\end{theorem}

Theorem~\ref{t2.5} will be established in Section~\ref{s3}. Here we
show that the exact controllability in a projection is a simple
consequence of UAC; in the next subsection, we deduce some corollaries
from Theorems~\ref{t2.2} and~\ref{t2.5}.

\begin{proof}[Proof of Theorem~\ref{t2.2}]
Let us fix a time~$T>0$, an initial point~$u_0\in V$, and a projection
${\mathsf P}_F:H\to H$ onto a finite-dimensional 
subspace~$F\subset H$.  Recall that~$B_F(R)$ stands for the closed
ball in~$F$ of radius~$R$ centred at origin and denote by~$C$ the norm
of~${\mathsf P}_F:H\to H$. Let us fix any~$R>C$ and choose~$\delta>0$
so small that
\begin{equation} \label{2.06}
\sup_{\hat u\in B_F(R)}\|e^{-\delta L}\hat u-\hat u\|\le\tfrac12.
\end{equation}
Denote by~$\KK$ the image of~$B_F(R)$ under~$e^{-\delta L}$. This is a
compact subset of~$V$, and by Theorem~\ref{t2.5}, there
is a continuous mapping\footnote{Recall that the control space
$C^\infty(J_T,E)\cap\Theta_T(h,u_0)$ is endowed with the metric
generated by the norm in~$L^2(J_T,H)$.}
$$
\varPsi:\KK\to C^\infty(J_T,E)\cap\Theta_T(h,u_0)
$$
such that 
\begin{equation} \label{2.07}
\sup_{\hat v\in\KK}\|\RR_T(u_0,\varPsi(\hat v))-\hat v\|_V<\tfrac12.
\end{equation}
It follows~\eqref{2.06} and~\eqref{2.07} that 
$$
\sup_{\hat u\in B_F(R)}
\|\RR_T(u_0,\varPsi(e^{-\delta L}\hat u))-\hat u\|<1.
$$
Therefore the continuous mapping
$$
\varPhi:B_F(R)\to F, \quad 
\hat u\mapsto{\mathsf P}_F\RR_T(u_0,\varPsi(e^{-\delta L}\hat u)),
$$
satisfies inequality~\eqref{1.1} with $\e=C$. Hence, by
Proposition~\ref{p1.1}, we have $\varPhi(B_F(R))\supset B_F(R-C)$.  In
particular, it follows that for any $\hat u\in B_F(R-C)$ there is
$\eta\in C^\infty(J_T,E)\cap\Theta_T(h,u_0)$ such that 
${\mathsf P}_F\RR_T(u_0,\eta)=\hat u$. Since $R>C$ is arbitrary, we
obtain the conclusion of Theorem~\ref{t2.2}.
\end{proof}

\subsection{\sloppy Solid controllability and Cauchy problem for the 
NS system}
\label{s2.3}
In this subsection, we establish some corollaries of
Theorems~\ref{t2.2} and~\ref{t2.5}. Let $E\subset U$ and $F\subset H$
be finite-dimensional subspaces, let~${\mathsf P}_F:H\to H$ be a
projection onto~$F$, and let~$T>0$ be a constant. 

\begin{definition} \label{d2.6}
The control system~\eqref{2.1} with $\eta\in L^2(J_T,E)$ is said to be
{\it solidly ${\mathsf P}_F$-controllable in time~$T$\/} if for any~$R>0$
and $u_0\in V$ there is a constant~$\e>0$ and a compact set
$\CC\subset L^2(J_T,E)\cap\Theta_T(h,u_0)$ such that, for any
continuous mapping $S:\CC\to F$ satisfying the inequality
\begin{equation} \label{2.7}
\sup_{\eta\in\CC}\|S(\eta)-{\mathsf P}_F\RR_T(u_0,\eta)\|_F\le\e,
\end{equation}
we have $S(\CC)\supset B_F(R)$. 
\end{definition}

\begin{proposition} \label{p2.7}
Under the conditions of Theorem~\ref{t2.2}, for any~$T>0$, any
finite-dimensional subspace $F\subset H$, and any projection
${\mathsf P}_F:H\to H$ onto~$F$, Eq.~\eqref{2.1} with 
$\eta\in L^2(J_T,E)$ is solidly ${\mathsf P}_F$-controllable in time~$T>0$.
\end{proposition}

\begin{proof}
Let us fix any constant~$R>0$, function~$u_0\in V$, and
subspace~$F\subset H$. As was shown in the proof of Theorem~\ref{t2.2}
(see Section~\ref{s2.2}), there is a continuous mapping
$\varPsi:B_F(R+2)\to L^2(J_T,E)\cap\Theta_T(h,u_0)$ such that
\begin{equation} \label{2.8}
\sup_{\hat u\in B_F(R+2)}\|{\mathsf P}_F\RR_T(u_0,\varPsi(\hat u))-\hat u\|
\le1.
\end{equation}
Let us set $\CC=\varPsi(B_F(R+2))$. Since $\dim F<\infty$
and~$\varPsi$ is continuous, we conclude that~$\CC$ is a compact
subset of $L^2(J_T,E)\cap\Theta_T(h,u_0)$.  Let $S:\CC\to F$ be an
arbitrary continuous mapping such that~\eqref{2.7} holds with
$\e=1$. Then it follows from~\eqref{2.8} that the mapping
$S\circ\varPsi:B_F(R+1)\to F$ satisfies the inequality
$$
\sup_{\hat u\in B_F(R+2)}\|S\circ\varPsi(\hat u)-\hat u\|_V\le2. 
$$
Applying Proposition~\ref{p1.1}, we see that
$S\circ\varPsi(B_F(R+2))\supset B_F(R)$.  It follows that
$S(\CC)\supset B_F(R)$. Since~$R>0$ was arbitrary, this completes the
proof of Proposition~\ref{p2.7}.
\end{proof}

We now show that the control function~$\eta$ in Theorem~\ref{t2.2} can
be taken from a finite-dimensional subspace. Namely, we have the
following result.

\begin{proposition} \label{p2.8}
Suppose that the conditions of Theorem~\ref{t2.2} are fulfilled, 
and let~$X$ be a vector space dense in~$L^2(J_T,E)$. Then for any positive
constants~$T$ and~$R$, any initial function~$u_0\in V$, any subspace
$F\subset H$ with $\dim F<\infty$, and any projection 
${\mathsf P}_F:H\to H$ onto~$F$, there is a ball~$B$ in a
finite-dimensional subspace $Y\subset X$ such that
\begin{equation} \label{2.9}
{\mathsf P}_F\RR_T(u_0,B)\supset B_F(R).
\end{equation}
In particular, we can take $X=C^\infty(J_T,E)$.
\end{proposition}

\begin{proof}
By Proposition~\ref{p2.7}, Eq.~\eqref{2.1} with $\eta\in L^2(J_T,E)$
is solidly ${\mathsf P}_F$-controllable in time~$T$.  Let $\e>0$ and
$\CC\subset L^2(J_T,E)\cap\Theta_T(h,u_0)$ be the corresponding
constant and compact set entering Definition~\ref{d2.6}.  It follows
from Proposition~\ref{p1.2} with $E=\{0\}$ and $v\equiv0$ that
$\Theta_T(h,u_0)$ is an open subset of~$L^2(J_T,E)$. Therefore there
is~$\delta>0$ such that
$$
O_\delta(\CC)=\bigl\{\eta\in L^2(J_T,E):\dist(\eta,\CC)\le\delta\bigr\}
\subset\Theta_T(h,u_0),
$$
where we set
$$
\dist(\eta,\CC)=\inf_{\zeta\in\CC}\|\eta-\zeta\|_{L^2(J_T,H)}.
$$
Furthermore, since~$X$ is dense in~$L^2(J_T,E)$, we can find a 
finite-dimensional subspace $Y\subset X$ such that
\begin{equation} \label{2.10}
\sup_{\eta\in\CC}\|P_Y\eta-\eta\|_{L^2(J_T,E)}\le\delta,
\end{equation}
where~$P_Y$ denotes the orthogonal projection in~$L^2(J_T,E)$ onto~$Y$.
It follows that
\begin{equation} \label{2.11}
P_Y\CC\subset O_\delta(\CC)\subset \Theta_T(h,u_0).
\end{equation}
By Proposition~\ref{p1.2}, the operator
$\RR(u_0,\cdot):\Theta_T(h,u_0)\to\XX_T$ is locally Lipschitz
continuous. Therefore, taking~$\delta>0$ sufficiently small, we deduce
from~\eqref{2.10} and~\eqref{2.11} that
$$
\sup_{\eta\in\CC}\|\RR_T(u_0,P_Y\eta)-\RR_T(u_0,\eta)\|_V\le\frac{\e}{C},
$$
where~$C$ is the norm of~${\mathsf P}_F$ regarded as an operator from~$V$ to~$H$.
Thus, the mapping $S(\eta)={\mathsf P}_F\RR_T(u_0,P_Y\eta)$ satisfies~\eqref{2.7}.
Hence, by Proposition~\ref{p2.7}, we have ${\mathsf P}_F\RR_T(u_0,P_Y\CC)\supset B_F(R)$.
It remains to note that~$P_Y\CC$ is contained in a ball of the finite-dimensional 
space~$Y\subset X$.
\end{proof}

We now consider the Cauchy problem for the NS equation~\eqref{1.2}.
Let $G\subset H$ be a closed vector space.  For any $u_0\in V$, $T>0$,
and~$\nu>0$, let~$\Xi_{T,\nu}(G,u_0)$ be the set of functions $f\in
L^2(J_T,G)$ for which problem~\eqref{1.2}, \eqref{2.2} has a unique
solution $u\in\XX_T$. If~$E\subset G$ is a closed subspace, then we
denote by~$G\ominus E$ the orthogonal complement of~$E$ in~$G$ and
by~$Q(T,G,E)$ the orthogonal projection in~$L^2(J_T,G)$ onto the
subspace~$L^2(J_T,G\ominus E)$.

\begin{proposition} \label{p2.9}
Let $E\subset U$ be a finite-dimensional subspace such that~$E_\infty$
is dense in~$H$ and let~$G\subset H$ be a closed subspace
containing~$E$.  Then~$\Xi_{T,\nu}(G,u_0)$ is a non-empty open subset
of~$L^2(J_T,G)$ such that
$$
Q(T,G,E)\Xi_{T,\nu}(G,u_0)=L^2(J_T,G\ominus E)\quad
\mbox{for any $T>0$, $\nu>0$, $u_0\in V$}.
$$
\end{proposition}

\begin{proof}
The fact that $\Xi_{T,\nu}(G,u_0)$ is open follows immediately from
Proposition~\ref{p1.2}.  The other claims of the proposition are
equivalent to the following property: 
for any $h\in L^2(J_T,G\ominus E)$ there is $\eta\in L^2(J_T,E)$ such
that $h+\eta\in\Xi_{T,\nu}(G,u_0)$. This is a straightforward
consequence of Theorem~\ref{t2.2}.
\end{proof}

\section{Proof of Theorem~\ref{t2.5}}
\label{s3}
\subsection{Scheme of the proof}
\label{s3.1}
Let~$E$ be a finite-dimensional vector space and let $E_1=\FF(E)$
(see~\eqref{2.5}).  Along with Eq.~\eqref{2.1}, consider two other
control systems:
\begin{align}
\dot u+\nu L(u+\zeta(t))+B(u+\zeta(t))&=h(t)+\eta(t),\label{3.1}\\
 \dot u+\nu Lu+B(u)&=h(t)+\eta_1(t). \label{3.2}
\end{align}
Here~$\eta$ and~$\zeta$ are $E$-valued controls and~$\eta_1$ is an
$E_1$-valued control. Let us fix a constant~$\e>0$, an initial point
$u_0\in V$, a compact set~$\KK\subset V$, and a vector space 
$X\subset L^2(J,H)$.  Equation~\eqref{2.1} with $\eta\in X$ is said to
be {\it uniformly $(\e,u_0,\KK)$-controllable\/} if there is a
continuous mapping
$$
\varPsi:\KK\to X\cap\Theta_T(h,u_0)
$$
such that~\eqref{2.6} holds. In what follows, if~$\e$, $u_0$,
and~$\KK$ are fixed in advance, then the above property will be called
{\it uniform $\e$-controllability\/}.

The concept of uniform $\e$-controllability for~\eqref{3.1} is defined
in a similar way. Namely, let~$\widehat\Theta_T(h,u_0)$ be the
set of pairs $(\eta,\zeta)\in L^2(J_T,H)\times L^4(J_T,H^2)$ for which
problem~\eqref{3.1}, \eqref{2.2} has a unique solution $u\in\XX_T$ and
let~$\widehat\RR$ be an operator that is defined on the set 
$$
D(\widehat\RR)=\bigl\{
(u_0,\eta,\zeta)\in V\times L^2(J_T,H)\times L^4(J_T,H^2): 
(\eta,\zeta)\in\widehat\Theta_T(h,u_0)\bigr\}
$$
and takes each triple~$(u_0,\eta,\zeta)\in D(\widehat\RR)$
to the solution~$u\in\XX_T$ of~\eqref{3.1}, \eqref{2.2}. Rewriting
Eq.~\eqref{3.1} in the form
$$
\dot u+\nu Lu+B(u)+B(u,\zeta)+B(\zeta,u)
=h(t)+\eta-\nu L\zeta-B(\zeta)
$$
and applying Proposition~\ref{p1.2} with $E=\{0\}$, we see that
$D(\widehat\RR)$ is an open subset of $V\times L^2(J_T,H)\times
L^4(J_T,H^2)$, and the operator~$\widehat\RR$ is locally Lipschitz
continuous on~$D(\widehat\RR)$. 

Now let $\widehat X\subset L^2(J,H)\times L^4(J,H^2)$ be a vector 
space, not necessarily closed. Equation~\eqref{3.1} with 
$(\eta,\zeta)\in\widehat X$ is said to be {\it uniformly 
$(\e,u_0,\KK)$-controllable\/} if there is a continuous mapping
$$
\widehat\varPsi:\KK\to \widehat X\cap\widehat\Theta_T(h,u_0)
$$
such that
\begin{equation} \label{3.3}
\|\widehat\RR_T(u_0,\widehat\varPsi(\hat u))-\hat u\|_V<\e\quad
\mbox{for any $\hat u\in\KK$},
\end{equation}
where~$\widehat\RR_t(u_0,\eta,\zeta)$ denotes the restriction 
of~$\widehat\RR(u_0,\eta,\zeta)$ to the time~$t$.

The proof of Theorem~\ref{t2.5} is based on the following three
propositions (cf.\ Propositions~3.1 and~3.2 and Section~2.2
in~\cite{S05}).  Let us fix a constant~$\e>0$, an initial
point~$u_0\in V$, and a compact subset~$\KK\subset V$.

\begin{proposition} \label{p3.1}
{\bf(extension principle)} 
Let~$E\subset U$ be a finite-dimensional vector space. Then
Eq.~\eqref{2.1} with $\eta\in C^\infty(J_T,E)$ is uniformly
$\e$-controllable if and only if so is Eq.~\eqref{3.1} with
$(\eta,\zeta)\in C^\infty(J_T,E\times E)$.
\end{proposition}

\begin{proposition} \label{p3.2}
{\bf(convexification principle)} Let~$E\subset U$ be a
finite-dimen\-sional subspace and let $E_1=\FF(E)$. Then
Eq.~\eqref{3.1} with $(\eta,\zeta)\in C^\infty(J_T,E\times E)$ is
uniformly $\e$-controllable if and only if so is Eq.~\eqref{3.2} with
$\eta_1\in C^\infty(J_T,E_1)$.
\end{proposition}

\begin{proposition} \label{p3.3}
Let~$E\subset U$ be a finite-dimensional vector space such
that~$E_\infty$ is dense in~$H$. Then there is an integer~$k\ge1$
depending on~$\e$, $u_0$, and~$\KK$ such that Eq.~\eqref{2.1} with
$\eta\in C^\infty(J_T,E_k)$ is uniformly $\e$-controllable.
\end{proposition}

If Propositions~\ref{p3.1}~--~\ref{p3.3} are established, then for
any~$\e>0$, $u_0\in V$, and~$\KK\subset V$ we first use
Proposition~\ref{p3.3} to find an integer~$k\ge1$ such that
Eq.~\eqref{2.1} with $\eta\in C^\infty(J_T,E_k)$ is uniformly
$(\e,u_0,\KK)$-controllable. Combining this property with
Propositions~\ref{p3.2} and~\ref{p3.1} in which $E=E_{k-1}$, we
conclude that Eq.~\eqref{2.1} with $\eta\in C^\infty(J_T,E_{k-1})$ is
uniformly $(\e,u_0,\KK)$-controllable.  Repeating this argument~$k-1$
times, we see that the same property is true for Eq.~\eqref{2.1} with
$\eta\in C^\infty(J_T,E)$. Since~$\e$, $u_0$, and~$\KK$ are arbitrary,
this completes the proof of Theorem~\ref{t2.5}.

To prove the above propositions, we repeat the scheme used
in~\cite{S05} (see Sections~2.2, 3.2, and~3.3). The important point
now is that we have to follow carefully the dependence of controls on
the final state~$\hat u$.  The proofs of
Propositions~\ref{p3.1}~--~\ref{p3.3} are carried out in next three
subsections. Here we formulate a lemma on uniform
$\e$-controllability; it will be used in
Sections~\ref{s3.2}~--~\ref{s3.4}. As before, we fix a
constant~$\e>0$, an initial point $u_0\in V$, and a compact
set~$\KK\subset V$.

\begin{lemma} \label{l3.4}
Let $X,Y\subset L^2(J_T,H)$ be vector spaces such that~$X$ is
contained in the closure of~$Y$ and Eq.~\eqref{2.1} with~$\eta\in X$
is uniformly $\e$-controllable. Then there is a finite-dimensional
subspace~$Y_0\subset Y$ such that Eq.~\eqref{2.1} with~$\eta\in Y_0$
is uniformly $\e$-controllable.
\end{lemma}

To prove this lemma, it suffices to repeat the argument used in the proof
of Proposition~\ref{p2.8}; we shall not dwell on it. Also note that an
analogue of Lemma~\ref{l3.4} is true for Eq.~\eqref{3.1}.

\subsection{Extension principle: proof of Proposition~\ref{p3.1}}
\label{s3.2}
We need to show that if Eq.~\eqref{3.1} with 
$(\eta,\zeta)\in C^\infty(J_T,E\times E)$ is uniformly
$\e$-controllable, then so is Eq.~\eqref{2.1} with 
$\eta\in C^\infty(J_T,E)$. Since $C^\infty(J_T,E)$ is dense
in~$L^2(J_T,E)$, in view of Lemma~\ref{l3.4}, it suffices to establish
that property for Eq.~\eqref{2.1} with $\eta\in L^2(J_T,E)$.

Recall that~$\mathsf P$ and~$\mathsf Q$ stand for the orthogonal projection
in~$H$ onto the subspaces~$E$ and~$E^\bot$, respectively. Let 
$$
\widehat\varPsi:\KK\to C^\infty(J_T,E\times E)\cap\widehat\Theta_T(h,u_0),
\quad \widehat\varPsi(\hat u)=\bigl(\eta(t,\hat u),\zeta(t,\hat u)\bigr),
$$
be an operator for which~\eqref{3.3} holds. We choose any sequence of functions
$\varphi_k\in C^\infty(\R)$ with the following properties:
\begin{gather} 
0\le\varphi_k(t)\le 1\quad\mbox{for all $t\in\R$},\label{3.4}\\
\varphi_k(t)=0\quad\mbox{for $t\le 0$ and $t\ge T$},\label{3.5}\\
\varphi_k(t)=1\quad\mbox{for $\tfrac1k\le t\le T-\frac1k$}.\label{3.6}
\end{gather}
We now define a sequence of continuous mappings 
$\varPsi_k:\KK\to L^2(J_T,E)$ by the following rule:
\begin{itemize}
\item for any $\hat u\in\KK$ and $k\ge1$, set
\begin{equation} \label{3.7}
v_k(t,\hat u)=\varphi_k(t)\zeta(t,\hat u)
+{\mathsf P}\widehat\RR_t(u_0,\widehat\varPsi(\hat u)), \quad t\in J_T;
\end{equation}
\item
denote by $w_k(\cdot,\hat u)\in\XX_T(E)$ the solution of the 
problem\,\footnote{We shall show that such a solution exists for $k\gg1$.}
\begin{equation} \label{3.8}
\begin{aligned}
\dot w+\nu L_Ew+{\mathsf Q}\bigl(B(w)+B(v_k,w)+B(w,v_k)\bigr)
&=f_k(t,\hat u),\\
w(0)&={\mathsf Q}u_0,
\end{aligned}
\end{equation}
where $v_k=v_k(t,\hat u)$ and 
\begin{equation} \label{3.9}
f_k(t,\hat u)={\mathsf Q}\bigl(h(t)-B(v_k(t,\hat u))-\nu Lv_k(t,\hat u)\bigr);
\end{equation}
\item
denote $u_k(t,\hat u)=w_k(t,\hat u)+w_k(t,\hat u)$ and define~$\varPsi_k$ 
by the formula
\begin{equation} \label{3.10}
\varPsi_k(\hat u)=\eta_k(t,\hat u)
:=\dot v_k+{\mathsf P}(\nu Lu_k+B(u_k)-h).
\end{equation}
\end{itemize}
We claim that for sufficiently large~$k\ge1$ the function~$\varPsi_k$
is well defined and continuous and satisfies~\eqref{2.6}. Indeed, let us write
$$
v(t,\hat u)={\mathsf P}\widehat\RR_t(u_0,\widehat\varPsi(\hat u)),\quad
w(t,\hat u)={\mathsf Q}\widehat\RR_t(u_0,\widehat\varPsi(\hat u)).
$$
Then $w(\cdot,\hat u)\in\XX_T(E)$ is a solution of the problem
\begin{equation} \label{3.11}
\begin{aligned}
\dot w+\nu L_Ew+{\mathsf Q}\bigl(B(w)+B(v,w)+B(w,v)\bigr)&=f(t,\hat u),\\ 
w(0)&={\mathsf Q}u_0,
\end{aligned}
\end{equation}
where $v=v(t,\hat u)$ and 
$f(t,\hat u)={\mathsf Q}\bigl(h(t)-B(v(t,\hat u))-\nu Lv(t,\hat u)\bigr)$.
We wish to consider~\eqref{3.8} as a perturbation of~\eqref{3.11}. 

Since~$\KK$ is compact and $\widehat\RR(u_0,\widehat\varPsi(\cdot)):\KK\to\XX_T$ 
is continuous, we have
\begin{equation} \label{3.12}
\sup_{\hat u\in\KK}
\bigl(\|v(\cdot,\hat u)\|_{\XX_T}+\|w(\cdot,\hat u)\|_{\XX_T}\bigr)<\infty.
\end{equation}
It follows from~\eqref{3.4}, \eqref{3.6}, and~\eqref{3.7} that
\begin{equation} \label{3.13}
\sup_{\hat u\in\KK}\|v_k(\cdot,\hat u)-v(\cdot,\hat u)\|_{L^4(J_T,U)}\to0
\quad\mbox{as $k\to\infty$}.
\end{equation}
Combining this with standard estimates for the nonlinear term and the fact
that $\dim E<\infty$, we conclude that (cf.~(3.8) in~\cite{S05})
\begin{equation} \label{3.14}
\sup_{\hat u\in\KK}\|f_k(\cdot,\hat u)-f(\cdot,\hat u)\|_{L^2(J_T,H)}\to0
\quad\mbox{as $k\to\infty$}.
\end{equation}
Proposition~\ref{p1.2} and relations~\eqref{3.12}~--~\eqref{3.14} imply that there 
is an integer~$k_0\ge1$ such that, for any $k\ge k_0$ and $\hat u\in\KK$,
problem~\eqref{3.8} has a unique solution $w_k(\cdot,\hat u)\in\XX_T(E)$. Moreover, 
the function $\hat u\mapsto w_k(\cdot,\hat u)$ is continuous from~$\KK$ to~$\XX_T(E)$.
It follows from~\eqref{3.10} that the operator~$\varPsi_k$ is well defined and 
continuous for~$k\ge k_0$. 

Let us show that~$\varPsi_k$ satisfies~\eqref{2.6} for $k\gg1$. Since
the resolving operator associated with~\eqref{3.11} is locally
Lipschitz continuous (see Proposition~\ref{p1.2}), 
for any $\hat u\in\KK$ and $k\ge k_0$, we have
$$
\|w_k(\cdot,\hat u)-w(\cdot,\hat u)\|_{\XX_T}
\le C\bigl(\|f_k(\cdot,\hat u)-f(\cdot,\hat u)\|_{L^2(J_T,H)}
+\|v_k(\cdot,\hat u)-v(\cdot,\hat u)\|_{L^4(J_T,V)}\bigr),
$$
where $C>0$ does not depend on~$k$ and~$\hat u$. Combining this inequality 
with~\eqref{3.13} and~\eqref{3.14}, we derive
\begin{equation} \label{3.15}
\sup_{\hat u\in\KK}\|w_k(\cdot,\hat u)-w(\cdot,\hat u)\|_{\XX_T}
\to0\quad\mbox{as $k\to\infty$}.
\end{equation}
Now note that, in view of~\eqref{3.5} and~\eqref{3.7}, we have
\begin{align*}
\|\RR_T(u_0,\varPsi_k(\hat u))-\hat u\|_V
&\le \|\RR_T(u_0,\varPsi_k(\hat u))-\widehat\RR_T(u_0,\widehat\varPsi_k(\hat u))\|_V\\
&\qquad+\|\widehat\RR_T(u_0,\widehat\varPsi_k(\hat u))-\hat u\|_V\\
&\le \|w_k(\cdot,\hat u)-w(\cdot,\hat u)\|_V
+\|\widehat\RR_T(u_0,\widehat\varPsi_k(\hat u))-\hat u\|_V.
\end{align*}
Taking the supremum over $\hat u\in\KK$ and using~\eqref{3.3} 
and~\eqref{3.15}, we see that~$\varPsi_k$ satisfies~\eqref{2.6} for 
sufficiently large~$k\ge k_0$. The proof of Proposition~\ref{p3.1}
is complete. 

\subsection{Convexification principle: proof of Proposition~\ref{p3.2}}
\label{s3.3}
We need to prove that if Eq.~\eqref{3.2} with $\eta_1\in C^\infty(J_T,E_1)$
is uniformly $\e$-con\-trollable, then so is Eq.~\eqref{3.1} with 
$(\eta,\zeta)\in L^2(J_T,E)\times L^4(J_T,E)$; the converse assertion is
obvious in view of Proposition~\ref{p3.1}. Let us outline the main idea. 

Let $\varPsi_1:\KK\to L^2(J_T,E_1)$ be a continuous mapping such that 
\begin{equation} \label{3.17}
\hat\e:=\sup_{\hat u\in\KK}\|\RR_T(u_0,\varPsi_1(\hat u))-\hat u\|_V<\e.
\end{equation}
By definition, the function 
$u_1(t,\hat u)=\RR_t(u_0,\varPsi_1(\hat u))$ satisfies Eq.~\eqref{3.2}
in which $\eta_1=\eta_1(\cdot,\hat u):=\varPsi_1(\hat u)$.  We wish to
approximate~$u_1(t,\hat u)$ by a solution~$u(t,\hat u)$ of
Eq.~\eqref{3.1} with some functions $\eta(\cdot,\hat
u),\zeta(\cdot,\hat u)\in L^\infty(J_T,E)$. This approximation should
be such that
\begin{equation} \label{3.18}
\sup_{\hat u\in\KK}\|u(T,\hat u))-u_1(T,\hat u)\|_V\le\e-\hat\e,
\end{equation}
and the mapping $\hat u\mapsto(\eta(\cdot,\hat u),\zeta(\cdot,\hat u))$
is continuous as an operator from~$\KK$ to the space 
$L^2(J_T,E)\times L^4(J_T,E)$. 

To construct~$u(t,\hat u)$, one could try to apply the argument used
in Section~3.3 of~\cite{S05} for approximating individual
solutions. Unfortunately, it does not work because it is difficult to
ensure that the resulting control functions~$\eta$ and~$\zeta$
continuously depend on~$\hat u$. To overcome this difficulty, we first
approximate~$\eta_1(t,\hat u)$ by a family of piecewise constant
controls~$\tilde\eta_1(t,\hat u)$ with range in the convex envelope of
a finite set not depending on~$\hat u$ (cf.\ Section~12.3
in~\cite{AS-2004}). We next repeat the scheme of~\cite{S05} to
construct an approximation for solutions~$\tilde u_1(t,\hat u)$
corresponding to~$\tilde\eta_1(t,\hat u)$. A difficult point of the
proof is to follow the dependence of the control functions 
on~$\hat u$.  In what follows, we shall omit the tilde from the
notation.

\medskip
The realisaion of the above scheme is divided into several steps. We
begin with a generalisation of the concept of uniform approximate
controllability.

\smallskip
{\it Step~1}. Let~$A=\{\eta_1^l,l=1,\dots,m\}\subset E_1$ be a finite set. 
For any integer~$s\ge1$, denote by~$P_s(J_T,A)$ the set of functions 
$\eta_1\in L^2(J_T,E_1)$ satisfying the following properties:
\begin{itemize}
\item
there are non-negative functions $\varphi_l\in L^\infty(J_T)$, $l=1\dots,m$,
such that
$$
\sum_{l=1}^m\varphi_l(t)=1,\quad \eta_1(t)=\sum_{l=1}^m\varphi_l(t)\eta_1^l
\quad\mbox{for $0\le t<T$};
$$
\item
the functions~$\varphi_l$ are representable in the form
$$
\varphi_l(t)=\sum_{r=0}^{s-1}c_{lr}I_{r,s}(t)\quad\mbox{for $0\le t<T$},
$$
where $c_{lr}\ge0$ are some constants and~$I_{r,s}$ denotes the indicator function
of the interval $J_T(r,s)=[t_r,t_{r+1})$ with $t_r=rT/s$.
\end{itemize}
The set $P_s(J_T,A)$ is endowed with the metric 
$$
d_P(\eta_1,\zeta_1)=\sum_{l=1}^m\|\varphi_l-\psi_l\|_{L^\infty(J_T)},
\quad\eta_1,\zeta_1\in P_s(J_T,A),
$$
where $\{\varphi_l\}$ and $\{\psi_l\}$ are the families of functions corresponding to~$\eta_1$ 
and~$\zeta_1$, respectively. 

Recall that we have fixed a constant~$\e>0$, an initial 
point~$u_0\in V$, and a compact set $\KK\subset V$. We shall say that
Eq.~\eqref{3.2} with $\eta_1\in P_s(J_T,A)$ is {\it uniformly
$\e$-controllable\/} if there is a continuous\footnote{We emphasise
that, in contrast to Definition~\ref{d2.4} in which the vector
space~$X$ is endowed with the norm of~$L^2(J_T,H)$, the
mapping~$\varPsi_1$ is required to be continuous with respect to the
topology of~$P_s(J_T,A)$, which is stronger than that
of~$L^2(J_T,H)$.}  mapping $\varPsi_1:\KK\to P_s(J_T,A)$ such that
$\varPsi(\hat u)\in \Theta_T(h,u_0)$ for any $\hat u\in\KK$,
and~\eqref{3.17} holds. A proof of the following lemma is based on a
standard argument of the control theory and is given in the Appendix.

\begin{lemma} \label{l3.5}
Let us assume that Eq.~\eqref{3.2} with $\eta\in C^\infty(J_T,E_1)$  is uniformly
$\e$-controllable. Then there is a finite set $A=\{\eta_1^l,l=1,\dots,m\}\subset E_1$
and an integer~$s\ge1$ such that Eq.~\eqref{3.2} with $\eta\in P_s(J_T,A)$
is uniformly $\e$-controllable.
\end{lemma}

Let $\varPsi_1:\KK\to P_s(J_T,A)$ be the function constructed in Lemma~\ref{l3.5}. 
We write
$$
\varPsi_1(\hat u)=\eta_1(t,\hat u)=\sum_{l=1}^m\varphi_l(t,\hat u)\eta_1^l.
$$
The definition of the space $P_s(J_T,A)$ and of its metric imply
that the functions~$\varphi_l$ have the form
\begin{equation} \label{3.21}
\varphi_l(t,\hat u)=\sum_{r=0}^{s-1}c_{lr}(\hat u)I_{r,s}(t),
\end{equation}
where $c_{lr}:\KK\to\R$ are non-negative continuous functions such that
$$
\sum_{l=1}^mc_{lr}(\hat u)=1\quad\mbox{for any $\hat u\in\KK$}.
$$
Since~$\eta_1^l\in\FF(E)$, by Lemma~3.3 in~\cite{S05}, there are vectors
$\eta^l,\zeta^{1l},\dots,\zeta^{kl}\in E$ and non-negative constants 
$\lambda_{1l},\dots,\lambda_{kl}$ such that $\sum_j\lambda_{jl}=1$ and 
$$
B(u)-\eta_1^l=\sum_{j=1}^k\lambda_{jl}
\bigl(B(u_1+\zeta^{jl})+\nu L\zeta^{jl}\bigr)-\eta^l\quad\mbox{for $u\in V$}.
$$
It follows that the function $u_1(\cdot)=\RR(u_0,\varPsi_1(\hat u))$ is a solution
of the equation
\begin{multline} \label{3.19}
\p_tu_1+\nu Lu_1+\sum_{j=1}^k\sum_{l=1}^m\lambda_{jl}\varphi_l(t,\hat u)
\bigl(B(u_1+\zeta^{jl})+\nu L\zeta^{jl}\bigr)\\
=h(t)+\sum_{j=1}^k\sum_{l=1}^m\lambda_{jl}\varphi_l(t,\hat u)\eta^l.
\end{multline}
Indexing the pairs $(j,l)$ by a single sequence $i=1,\dots,q$, 
we can write~\eqref{3.19} as
\begin{equation} \label{3.20}
\p_tu_1+\nu L\,\biggl(u_1+\sum_{i=1}^q\psi_i(t,\hat u)\zeta^i\biggr)
+\sum_{i=1}^q\psi_i(t,\hat u)B(u_1+\zeta^i)=h(t)+\eta(t,\hat u).
\end{equation}
Here $\zeta^i\in E$, $i=1,\dots,q$, are some vectors, 
$\eta(\cdot,\hat u)$ denotes the sum on the right-hand side
of~\eqref{3.19}, and
\begin{equation} \label{3.22}
\psi_i(t,\hat u)=\sum_{r=0}^{s-1}d_{ir}(\hat u)I_{r,s}(t),
\end{equation}
where $d_{ir}\in C(\KK)$ are non-negative functions such that 
$\sum_id_{ir}\equiv1$.

\smallskip
{\it Step~2}. We now approximate~$u_1(t,\cdot)$ by solutions of
Eq.~\eqref{3.1}.  To this end, we first assume that there is only one
interval of constancy, that is, $s=1$. In this case, the sums
in~\eqref{3.21} and~\eqref{3.22} contain only one term, and
Eq.~\eqref{3.20} takes the form
\begin{equation} \label{3.23}
\p_tu_1+\nu L\,\biggl(u_1+\sum_{i=1}^qd_i(\hat u)\zeta^i\biggr)
+\sum_{i=1}^qd_i(\hat u)B(u_1+\zeta^i)=h(t)+\eta(t,\hat u),
\end{equation}
where $d_i\in C(\KK)$ and $\eta\in C(\KK,L^2(J_T,E))$. Let us fix an
integer $k\ge1$ and, following a classical idea in the control theory,
define a sequence of continuous mappings $\zeta_k:\KK\to L^4(J_T,H^2)$
as
\begin{equation} \label{3.24}
\zeta_k(t,\hat u)=\zeta(kt/T,\hat u),
\end{equation}
where $\zeta(\cdot,\hat u)$ is a $1$-periodic function on~$\R$ such
that  
\begin{equation} \label{3.25}
\zeta(s,\hat u)=\zeta^i\quad
\mbox{for $0\le s-(d_1(\hat u)+\cdots+d_{i-1}(\hat u))<d_i(\hat u)$,
$i=1,\dots,q$}.
\end{equation}
It is easy to see that $\{\zeta_k(\cdot,\hat u), \hat u\in\KK,k\ge1\}$
is a bounded subset in $L^\infty(J_T,E)$. Let us rewrite~\eqref{3.23}
in the form
\begin{equation} \label{3.26}
\p_tu_1+\nu L(u_1+\zeta_k(t,\hat u))+B(u_1+\zeta_k(t,\hat u))
=h(t)+\eta(t,\hat u)+f_k(t,\hat u),
\end{equation}
where $f_k(t,\hat u)=f_{k1}(t,\hat u)+f_{k2}(t,\hat u)$, 
\begin{align} 
f_{k1}(t,\hat u)
&=\nu L\biggl(\zeta_k(t,\hat u)-\sum_{i=1}^qd_i(\hat u)\zeta^i\biggr), 
\label{3.27}\\
f_{k2}(t,\hat u)&=B(u_1(t,\hat u)+\zeta_k(t,\hat
u))-\sum_{i=1}^id_i(\hat u)B(u_1(t,\hat u)+\zeta^i).
\label{3.28}
\end{align}
We shall need the following result, which will be proved in the next
steps. Denote by~$B_V(u,r)$ the closed ball in~$V$ of radius~$r$
centred at~$u$. 

\begin{lemma} \label{l3.6}
For any $\e_0>0$ there is an integer~$k_0\ge1$ and a
constant~$\delta_0>0$ such that for any $k\ge k_0$, $\hat u\in\KK$,
and $v_0\in B_V(u_0,\delta_0)$, the problem
\begin{equation} \label{3.30}
\p_tv+\nu L(v+\zeta_k(t,\hat u))+B(v+\zeta_k(t,\hat u))
=h(t)+\eta(t,\hat u),\quad v(0)=v_0
\end{equation}
has a unique solution $v_k(\cdot,\hat u)\in\XX_T$, which satisfies the
inequality 
\begin{equation} \label{3.31}
\|v_k(\cdot,\hat u)-u_1(\cdot,\hat u)\|_{C(J_T,V)}\le\e_0.
\end{equation}
\end{lemma}

In particular, taking $\e_0=\hat\e$, where $\hat\e$ is the constant 
in~\eqref{3.17}, and defining the operator
$$
\widehat\varPsi_k:\KK\to L^2(J_T,E)\times L^4(J_T,E),\quad
\hat u\mapsto(\eta_k(\cdot,\hat u),\zeta_k(\cdot,\hat u)),
$$
we conclude that $\widehat\varPsi_k(\hat u)\in\widehat\Theta_T(h,v_0)$
for $v_0\in B_V(u_0,\e)$ and~$k\ge k_0$, and
$$
\sup_{\hat u\in\KK}
\|\widehat\RR_T(u_0,\widehat\varPsi_k(\hat u))-u_1(T,\hat u)\|_V
\le\hat\e\quad\mbox{for $k\ge k_0$}.
$$
Combining this with~\eqref{3.17}, we obtain
$$
\sup_{\hat u\in\KK}
\|\widehat\RR_T(u_0,\widehat\varPsi_k(\hat u))-\hat u\|_V
<\e\quad\mbox{for $k\ge k_0$}.
$$
Hence, Eq.~\eqref{3.1} with 
$(\eta,\zeta)\in L^2(J_T,E)\times L^4(J_T,E)$ is uniformly
$\e$-controllable.

\smallskip
{\it Step~3}. We now prove Lemma~\ref{l3.5}. Literal repetition of the
arguments in~\cite[Section~3.3]{S05} (see Step~2) shows that
the required assertion will be established if we prove the convergence
\begin{equation} \label{3.29}
\sup_{\hat u\in\KK}\bigl(\|K_\nu f_k(\cdot,\hat u)\|_{C(J_T,V)}
+\|B(K_\nu f_k(\cdot,\hat u))\|_{L^2(J_T,H)}\bigr)\to0
\quad\mbox{as $k\to\infty$},
\end{equation}
where
$$
K_\nu f(t)=\int_0^te^{-\nu(t-s)L}f(s)\,ds.
$$
Furthermore, in view of the calculations of Steps~3--4
in~\cite[Section~3.3]{S05}, it suffices to show that
\begin{equation} \label{3.32}
\sup_{\hat u\in\KK}\|F_k(\cdot,\hat u)\|_{C(J_T,H)}\to0
\quad\mbox{as $k\to\infty$},
\end{equation}
where $F_k(t,\hat u)=\int_0^tf_k(s,\hat u)\,ds$. To this end, we first
note that\footnote{See Step~5 in~\cite[Section~3.3]{S05}.} 
\begin{equation} \label{3.33}
\|F_k(\cdot,\hat u)\|_{C(J_T,H)}\to0\quad
\mbox{as $k\to\infty$ for any $\hat u\in\KK$}. 
\end{equation}
Suppose now that we have proved the uniform equicontinuity of the
family of mappings
\begin{equation} \label{3.34}
{\mathbf f}_k:\KK\to L^1(J_T,H),\quad 
\hat u\mapsto f_k(\cdot,\hat u).
\end{equation}
In this case, the family 
$\bigl\{\hat u\mapsto\int_0^\cdot f_k(s,\hat u)\,ds,k\ge1\bigr\}$ is
uniformly equicontinuous from~$\KK$ to~$C(J_T,V)$.  Combining this
property with~\eqref{3.33}, we arrive at~\eqref{3.32}.

\smallskip
{\it Step~4}.
We now show that~\eqref{3.34} is uniformly equicontinuous. The
explicit formulas~\eqref{3.27} and~\eqref{3.28} and standard estimates
for the bilinear form~$B$ show that it suffices to prove that the
function $\hat u\mapsto\zeta_k(\cdot,\hat u)$ is uniformly
equicontinuous from~$\KK$ to~$L^4(J_T,U)$. It follows from~\eqref{3.24}
and~\eqref{3.25} that
\begin{align*}
\|\zeta_k(\cdot,\hat u_1)-\zeta_k(\cdot,\hat u_2)\|_{L^4(J_T,U)}^4
&=\int_0^T\|\zeta(kt/T,\hat u_1)-\zeta(kt/T,\hat u_2)\|_U^4dt\\
&=T\int_0^1\|\zeta(s,\hat u_1)-\zeta(s,\hat u_2)\|_U^4ds\\
&\le C\sum_{i=1}^q|d_i(\hat u_1)-d_i(\hat u_2)|, 
\end{align*}
where $\hat u_1,\hat u_2\in\KK$ are arbitrary points and~$C>0$ is a
constant depending only on~$T$, $q$, and $\max_i\|\zeta^i\|_U$. Since
the functions~$d_i$ are uniformly continuous on the compact set~$\KK$,
we obtain the required result. This completes the proof of
Proposition~\ref{p3.2} in the case $s=1$.

\smallskip
{\it Step~5}.
We now consider the case of any $s\ge2$. Let us set
$I_r=[t_r,t_{r+1}]$ and $\XX^r=C(I_r,V)\cap L^2(I_r,U)$. For any
$r=0,\dots,s-1$, we denote by $\Theta^r(h,u_0)$ the set of functions
$(\eta,\zeta)\in L^2(I_r,H)\times L^4(I_r,H^2)$ for which
Eq.~\eqref{3.1} has a unique solution $u\in\XX^r$ satisfying the
initial condition
\begin{equation} \label{3.35}
u(t_r)=u_0.
\end{equation}
Introduce the set
$$
\DD^r=\bigl\{(u_0,\eta,\zeta)\in 
V\times L^2(I_r,H)\times L^4(I_r,H^2):
(\eta,\zeta)\in\Theta^r(h,u_0)\bigr\} 
$$
and define an operator $S^r:\DD^r\to V$ that takes each triple
$(u_0,\eta,\zeta)\in\DD^r$ to~$u(t_{r+1})$, where~$u\in\XX^r$ is the
solution of~\eqref{3.1}, \eqref{3.35}. It follows from
Proposition~\ref{p1.2} that the operator~$S^r$ is locally Lipschitz
continuous. 

We now define positive constants~$\beta_r$,
$r=0,\dots,s$, and continuous operators 
$\varPsi^r:\KK\to L^2(I_r,E)\times L^4(I_r,E)$, $r=0,\dots,s-1$, by
the following rule: 
\begin{itemize}
\item
set $\beta_s=\hat\e$, where~$\hat\e$ is the constant in~\eqref{3.17}; 
\item
if $\beta_{r+1}$ is constructed for some $r\le s-1$, then apply
Lemma~\ref{l3.6} with $\e_0=\beta_{r+1}$ to the interval~$I_r$ and
denote by~$\delta_0$ and~$k_0$ the corresponding parameters;
\item
set $\beta_r=\delta_0$ and $\varPsi^r=\widehat\varPsi_{k_0}$.
\end{itemize}
The construction implies that, for any $v_0\in B_V(u_1(t_r),\beta_r)$, 
$r=0,\dots,s-1$, and $\hat u\in\KK$, we have
\begin{equation} \label{3.36}
\varPsi^r(\hat u)\in\Theta^r(h,v_0), \quad
\|S^r(v_0,\varPsi^r(\hat u))-u_1(t_{r+1},\hat u)\|_V\le\beta_{r+1}.
\end{equation}
Let us define an operator 
$\widehat\varPsi:\KK\to L^2(J_T,E)\times L^4(J_T,E)$ as
$$
\widehat\varPsi(\hat u)(t)=\varPsi^r(\hat u)(t) \quad
\mbox{for $t\in I_r$, $r=0,\dots,s-1$}.
$$
It follows from~\eqref{3.36} that
$$
\widehat\varPsi(\hat u)\in\Theta_T(h,u_0),\quad
\|\RR_T(u_0,\widehat\varPsi(\hat u))-u_1(T,\hat u)\|_V\le\beta_s
=\hat\e\quad\mbox{for $\hat u\in\KK$}. 
$$
Comparing this with~\eqref{3.17}, we obtain~\eqref{2.6}. It remains to
note that since the functions~$\varPsi^r$ are continuous, so
is~$\widehat\varPsi$. This completes the proof of
Proposition~\ref{p3.2}.

\subsection{Proof of Proposition~\ref{p3.3}}
\label{s3.4}
{\it Step~1}. 
We first show that if the integer~$N\ge1$ is sufficiently
large, then Eq.~\eqref{2.1} with $\eta\in L^2(J_T,H_N)$ is uniformly 
$\e$-controllable. To this end, we fix a (small) constant~$\delta>0$ and define
a family of functions
\begin{equation} \label{3.50}
v_N(t,\hat u)
=T^{-1}{\mathsf P}_N\bigl(te^{-\delta L}\hat u+(T-t)e^{-tL}u_0\bigr),
\quad 0\le t\le T.
\end{equation}
It is easy to see that
\begin{align} 
K_\delta&:=\sup_{\hat u,N}\|v_N(\cdot,\hat u)\|_{\XX_T}<\infty
\quad\mbox{for any $\delta>0$},\label{3.51}\\
c_\delta&:=\sup_{\hat u,N}\|v_N(T,\hat u)-{\mathsf P}_N\hat u\|_V\to0
\quad\mbox{as $\delta\to0$},\label{3.52}
\end{align}
where the supremums are taken over $N\ge1$ and $\hat u\in\KK$.
We now choose a constant~$\delta>0$ so small that
\begin{equation} \label{3.53}
c_\delta\le\frac{\e}{3}.
\end{equation}

Consider the Cauchy problem
\begin{equation} \label{3.54}
\dot w+{\mathsf Q}_NL(w+v_N)+{\mathsf Q}_NB(w+v_N)={\mathsf Q}_Nh(t),
\quad w(0)={\mathsf Q}_Nu_0.
\end{equation}
Proposition~\ref{p1.4} and inequality~\eqref{3.51} imply
that there is an integer~$N_\delta\ge1$ not depending on~$\hat u\in\KK$
such that for any $N\ge N_\delta$ problem~\eqref{3.54} has a unique
solution $w_N(\cdot,\hat u)\in\XX_T(N_H)$. It follows that the function
$u_N(t,\hat u)=v_N+w_N$ belongs to~$\XX_T$ for any~$N\ge N_\delta$ and satisfies
Eq.~\eqref{2.1} with
\begin{equation} \label{3.49}
\eta(t)=\eta_N(t,\hat u):=\dot v_N+{\mathsf P}_N(Lu_N+B(u_N)-h).
\end{equation}
The required assertion will be established if we prove the following two claims:
\begin{itemize}
\item[\bf(a)]
For any~$N\ge N_\delta$, the function $\varPsi:\hat u\mapsto \eta_N(\cdot,\hat u)$
is continuous from~$\KK$ to~$L^2(J_T,H)$.
\item[\bf(b)]
We have
\begin{equation} \label{3.55}
\sup_{\hat u\in\KK}\|w_N(T,\hat u)\|_V\to0\quad\mbox{as $N\to\infty$}.
\end{equation}
\end{itemize}
Indeed, the very construction of~$\varPsi$ implies that
$$
\varPsi(\hat u)\in\Theta_T(h,u_0),\quad 
\RR(u_0,\varPsi(\hat u))=u_N. 
$$ 
Furthermore, it follows from~\eqref{3.52}, \eqref{3.53},
and~\eqref{3.55} that if~$N\ge N_\delta$ is sufficiently large, then
\begin{align} 
\sup_{\hat u\in\KK}\|\RR_T(u_0,\varPsi_N(\hat u))-\hat u\|_V
&\le c_\delta+\sup_{\hat u\in\KK}
\bigl(\|w_N(T,\hat u)\|_V+\|{\mathsf Q}_N\hat u\|_V\bigr)
\notag\\
&\le\frac{2\e}{3}+\sup_{\hat u\in\KK}\|{\mathsf Q}_N\hat u\|_V.
\label{3.56}
\end{align}
Since~$\KK\subset V$ is compact, the second term on the right-hand side
of~\eqref{3.56} can be made smaller than~$\frac\e3$ by choosing a sufficiently 
large~$N\ge N_\delta$.

\medskip
{\it Step~2}. Let us prove~(a) and~(b). 
Since~$\delta>0$, it follows from~\eqref{3.50} that the 
function $\hat u\mapsto v_N(\cdot,\hat u)$ is continuous from~$\KK$ 
to~$\XX_T$. By Proposition~\ref{p1.4}, the solution $w_N\in\XX_T(H_N)$ of
problem~\eqref{3.54} continuously depends on~$v_N\in\XX_T$. The continuity
of~$\varPsi$ follows now from~\eqref{3.49} and~\eqref{3.50}. 

The proof of~(b) literally repeats the argument used in~\cite{S05} (see the
proof of~(2.12)), and therefore we omit it.

\medskip
{\it Step~3}. 
We now show that if~$k\ge1$ is sufficiently large, then Eq.~\eqref{2.1}
with $\eta\in C^\infty(J_T,E_k)$ is uniformly $\e$-controllable. To this end,
we use Lemma~\ref{l3.4}.

Let us denote by~$Y\subset L^2(J_T,H)$ the union of the vector spaces 
$C^\infty(J_T,E_k)$, $k\ge1$. Since~$E_\infty$ is dense in~$H$, we conclude 
that~$L^2(J_T,H_N)$ is contained in the closure of~$Y$ for any~$N\ge1$. By
Lemma~\ref{l3.4}, there is a finite-dimensional subspace~$Y_0\subset Y$ such that
Eq.~\eqref{2.1} with~$\eta\in Y_0$ is uniformly $\e$-controllable. 
Since $\{C^\infty(J_T,E_k)\}_{k\ge1}$ is an increasing sequence of subspaces,
we see that $Y_0\subset C^\infty(J_T,E_k)$ for a sufficiently large~$k\ge1$.
This completes the proof of Proposition~\ref{p3.3}. 

\section{Appendix: proof of Lemma~\ref{l3.5}}
Let~$d$ be the dimension of~$E_1$ and let $\EE=\{e_1,\dots,e_d\}$ be a
basis in~$E_1$. We endow~$E_1$ with a scalar product~$(\cdot,\cdot)$
for which~$\EE$ is an orthonormal system. Let $\varPsi_1:\KK\to
C^\infty(J_T,E_1)$ be a continuous operator satisfying~\eqref{3.17}. In
view of Lemma~\ref{l3.4} (in which $X=Y=C^\infty(J_T,E_1)$), we can
assume without loss of generality that $\varPsi_1(\hat u)\in Y_0$ for
any $\hat u\in\KK$, where $Y_0\subset C^\infty(J_T,E_1)$ is a
finite-dimensional subspace.  
Let us set $\eta_1(\cdot,\hat u)=\varPsi_1(\hat u)$ and write
\begin{equation} \label{4.1}
\eta_1(t,\hat u)=\sum_{l=1}^d\zeta_l(t,\hat u)e_l,
\end{equation}
where $\zeta_l(t,\hat u)=(\eta_l(t,\hat u),e_l)$. Since all the norms
in the finite-dimensional space~$Y_0$ are equivalent, what has been
said implies that $\zeta_l\in C(J_T\times\KK)$ for $l=1,\dots,d$. Let
$$
M=\max_{l,t,\hat u}|\zeta_l(t,\hat u)|,
$$
where the maximum is taken over $l=1,\dots,d$ and 
$(t,\hat u)\in J_T\times\KK$.  We now set~$m=2d$,
$$
\eta_1^l=dMe_l\quad\mbox{for $l=1,\dots,d$}, \qquad
\eta_1^l=-dMe_l\quad\mbox{for $l=d+1,\dots,m$}. 
$$
In this case, we can rewrite~\eqref{4.1} in the form
$$
\eta_1(t,\hat u)=\sum_{l=1}^m\tilde\zeta_l(t,\hat u)\eta_1^l,
$$
where $\tilde\zeta_l\in C(J_T\times\KK)$, $l=1,\dots,m$, are
non-negative functions whose sum is equal to~$1$.

For any integer~$s\ge1$, let us set
$$
\varPsi_1^s(\hat u)=\sum_{l=1}^m\psi_{ls}(t,\hat u)\eta_1^l,
$$
where $\psi_{ls}(t,\hat u)=\tilde\zeta_l(rT/s,\hat u)$ for 
$t\in J_T(r,s)$. It is clear that $\varPsi_1^s(\cdot)$
is a continuous function from~$\KK$ to~$P_s(J_T,A)$, 
where $A=\{\eta_1^l,l=1,\dots,m\}$. Furthermore, since $\KK\subset V$ is compact, 
it is not difficult to show that
$$
\sup_{\hat u\in\KK}\|\varPsi_1^s(\hat u)-\varPsi_1(\hat u)\|_{L^2(J_T,H)}
\to0\quad\mbox{as $s\to\infty$}.
$$
Proposition~\ref{p1.2} now implies that $\varPsi_1^s(\hat u)\in\Theta_T(h,u_0)$ 
for any $\hat u\in\KK$ and sufficiently large~$s$, and we have
$$
\sup_{\hat u\in\KK}\|\RR_T(u_0,\varPsi_1^s(\hat u))-\RR_T(u_0,\varPsi_1(\hat u))\|_V
\to0\quad\mbox{as $s\to\infty$}.
$$
Combining this with~\eqref{3.17}, we conclude that Eq.~\eqref{3.2}
with $\eta_1\in P_s(J_T,A)$ is uniformly $\e$-controllable.

\end{document}